\documentclass[12pt,leqno,a4paper]{article}      

\usepackage{graphicx}

\usepackage{amsmath}
\usepackage{amsthm}
\usepackage{amssymb}
\usepackage{xypic}
\usepackage{mathrsfs} 

\theoremstyle{definition}
\newtheorem{exa}{Example}
\newtheorem{thm}[exa]{Theorem}
\newtheorem{lemma}[exa]{Lemma}
\newtheorem{prop}[exa]{Proposition}

\renewcommand{\P}{\mathbb P}   % For Projective space
\newcommand{\Q}{{\mathbb Q}}     % For rational numbers
\newcommand{\R}{{\mathbb R}}     % For real numbers
\newcommand{\C}{{\mathbb C}}     % For Complex numbers
\newcommand{\sss}{\scriptscriptstyle}
\newcommand{\tth}{\thinspace}

\newcommand{\dimost}{\noindent{\it Proof\ }}

\begin{document}

\title{\bf A remark on the Generalized Hodge 
Conjecture\tth{}\footnote{Supported by the research program 
PRIN \emph{Geometria delle variet\`a algebriche e dei loro spazi 
di moduli}, partially funded by MIUR (cofin 2008).}}

\author{\bf Dario Portelli}

\date{}
\maketitle

%%%%%%%%%%%%%%%%%%%%%%%%%%%%%%%%%%%%%%%%%%%%%%%%%%%%%%%%%%%%%%%%%%%%%%%%%%%%%%
%%%%%%%%%%%%%%%%%%%%%%%%%%%%%%%%%%%%%%%%%%%%%%%%%%%%%%%%%%%%%%%%%%%%%%%%%%%%%%

\section{Introduction}
\label{intro}

Let $X$ be a smooth, projective, connected algebraic variety over 
$\C ,$ of dimension $n.$   
We will simply write $H_r\tth X$ and $H^{i} X$ for the homology 
and cohomology groups with coefficients in $\Q\tth .$
The symbol \ $\subset$ \ will denote nonstrict inclusion.

\medskip

After Lefschetz proved his theorem on $(1,1)$-classes, one of the main 
problems of the trascendental theory of algebraic varieties was the 
characterization in $H_{{\sss 2}m}\tth X$ of the classes of  
algebraic subvarieties of $X,$ of (\tth pure\tth ) dimension $m\tth .$ 

Hodge proposed to face this problem from a broader point of view, and 
to this aim he introduced (\tth\cite{H}, p. 184\tth ) the following
subspaces of $H_r\tth X,$ which depend on the integer $p\geq 0$ 
\begin{equation}\label{preniveau}
\sum_{
Y\subset X \ closed,\tth {\bf algebraic} \atop 
{\rm codim}_{{}_X}Y\tth\geq\tth p
}
Im\bigl(\tth H_r\tth Y\to H_r\tth X\tth\bigr)
\end{equation}
The map $H_r\tth Y\to H_r\tth X$ here is the canonical one, 
induced by the inclusion $Y\subset X.$ These subspaces form 
a decreasing filtration on $H_r\tth X.$ In particular,
if $r=2m$ and $p=n-m\tth ,$ the corresponding space of the 
filtration consists of the classes of $m$-dimensional algebraic 
subvarieties of $X.$ 

Returning to the general case, Hodge found a necessary condition, in
terms of suitable integrals, for a homology 
$r$-class to be contained into the space (\ref{preniveau}), 
and asked if this condition is also sufficient. 
Grothendieck showed that the answer is negative \cite{G}. The main 
step for this was the translation of the whole set-up from homology
to cohomology. Let us briefly describe how the image of the
space (\ref{preniveau}) by the Poincar\'e duality map $PD:H_r\tth X\to
H^{\tth{\sss 2}n-r}\tth X$ can be defined directly in cohomology.

Consider a closed, algebraic subset $Y\subset X;$ to simplify,
we will assume $Y$ equidimensional, of codimension $p\tth .$ 
We will always consider algebraic subsets of $X$ with their reduced structure.
Therefore we can apply to $Y$ Hironaka's theorem on the 
resolution of singularities. If $\rho :\widetilde Y\to Y$ is a resolution 
of the singularities of $Y,$ let $f:\widetilde Y\to X$ denote the composition
of $\rho$ with the inclusion; $f$ is a proper map. 
The following diagram is trivially commutative because the Gysin map
$f_*$ in it is defined as the composition of the other maps (\tth throughout 
the paper we will use the convention \tth $i:=2n-r$\tth ) 
\begin{equation}\label{pesi}
\diagram
H^{\tth i-{\sss 2}p}\tth\widetilde Y\dto_{PD}\rrto^{f_*}&& H^i\tth X\\
H_r\tth\widetilde Y\rto_{\rho_*} & H_r\tth Y\rto & H_r\tth X\uto_{PD}
\enddiagram
\end{equation}
By the universal coefficient theorem \tth $H_r\tth X=Hom_{{}_\Q}(H^r\tth X,\Q )$
and similarly for $\widetilde Y$ and $Y.$ Therefore $H_r\tth\widetilde Y$ and 
$H_r\tth X$ are pure rational Hodge structures of weight $-r\tth ,$ and
$H_r\tth Y$ is a mixed rational Hodge structure. Then, a simple weight 
argument shows that
$$
Im\bigl(\tth H_r\tth\widetilde Y\to H_r\tth X\tth\bigr)\tth =\tth
Im\bigl(\tth H_r\tth Y\to H_r\tth X\tth\bigr)
$$
More generally, consider a map $f:Z\to X,$ with $Z$ projective, smooth,
connected and \tth ${\rm dim}\tth Z\leq n-p\tth .$ Then we have a commutative
diagram like (\ref{pesi}), with $Z$ instead of $\widetilde Y$ and $f(Z)$
instead of $Y,$ and the same weight argument shows that
$$
Im(f_*)\tth =\tth PD\bigl(\tth Im\bigl(\tth H_r\tth f(Z)\to 
H_r\tth X\tth\bigr)\tth\bigr)
$$
Therefore, the image by Poincar\'e duality of the space (\ref{preniveau})   
is the space $N^{\tth p}H^{\tth i}\tth X$ of the coniveau filtration of 
$H^{\tth i}\tth X\tth .$ To summarize, a class 
$\xi\in H^{\tth i}\tth X$ is in $N^{\tth p}H^{\tth i}\tth X$
if and only if there are a closed algebraic subset $Y\subset X,$ 
of codimension $\geq p\tth ,$ and a $(2n-i)$-cycle $\Gamma\subset
Y$ (\tth this forces \ $p\leq i/2$\tth ) such that $PD(\xi)=[\Gamma ]\tth .$  

\medskip

Let $Z$ be as above, with \tth ${\rm dim}\tth Z=m\leq n-p\tth .$ In the up-left corner 
of (\ref{pesi}) we have $H^{\tth{\sss 2}m-r}\tth Z;$ since \tth
$2m-r\leq 2(n-p)-r=i-2p\tth ,$ we have \tth $Im(f_*)\subset F^{\tth p}
H^{\tth i}(X,\C )$ because the Gysin map $f_*$ is a morphism of Hodge structures. 
Hence
\begin{equation}\label{GHC0}
N^{\tth p}H^{\tth i}\tth X\tth\subset\tth
F^{\tth p}H^{\tth i}(X,\C )\tth\cap\tth H^{\tth i}\tth X
\end{equation}
This corresponds to the translation in cohomology of Hodge necessary condition 
for homology $r$-classes to be contained into the space (\ref{preniveau})
mentioned above. 
Hodge's original problem, solved by Grothendieck for the negative, was whether 
this inclusion is an equality or not. Over the years this problem became known 
as the Generalized Hodge Conjecture. The aim of this paper is to remark that it 
is possible to replace $F^{\tth p}H^{\tth i}(X,\C )$ in (\ref{GHC0}) by a canonically 
defined subspace \tth $S^{\tth p,i}$ of $F^{\tth p}H^{\tth i}(X,\C )\tth ,$ 
such that the refined (\ref{GHC0}) obtained in this way is always an equality.

\medskip  

Before defining the $S^{\tth p,i},$ let me anticipate 
what follows. We will always denote by $\Delta_r\subset\R^r$ the standard, open 
$r$-simplex, and by $\bar\Delta_r$ its closure. In the next section we will 
prove in Lemma \ref{omologia} that every class in $H_r\tth X$ can be represented 
by a singular $r$-cycle $\Gamma = \sum_h\tth m_h\tth\sigma_h$ such that the
restriction to $\Delta_r$ of every singular simplex $\sigma_h:\bar\Delta_r\to X$ 
appearing in $\Gamma$ is a real-analytic embedding. It will be clear from the
proof of Lemma \ref{omologia} that the same is true for the classes of
$H_r\tth Y,$ for any closed, algebraic set $Y\subset X.$
When dealing with singular
simplexes we will {\sl always} refer to their restriction to the interior $\Delta_r$
of $\bar\Delta_r\tth .$ 

\medskip 

Fix a class $\xi\in N^{\tth p}H^{\tth i}\tth X$ and let
$Y\subset X$ be a closed algebraic set, of codimension $\geq p\tth ,$ 
and $\Gamma\subset Y$ a $r$-cycle such that $PD(\xi)=[\Gamma ]\tth .$ 
Fix a singular simplex $\sigma_h:\bar\Delta_r\to X$ of $\Gamma ,$ and denote
$\sigma_h(\Delta_r)$ by $\Delta\tth .$ If $\Delta$ is contained in the
singular locus $Y_{sing}$ of $Y,$ then replace $Y$ by $Y_{sing}\tth .$
Decompose $Y_{sing}$ into its non-singular and singular parts.  
If $\Delta$ is contained into $(Y_{sing})_{sing}\tth ,$ then replace
$Y_{sing}$ by $(Y_{sing})_{sing}\tth .$
Iterating this procedure, we can always reduce
to the case when $\Delta$ is contained in an algebraic subset $Y'$ of $X,$
of codimension $\geq p\tth ,$ and $\Delta$ is not contained into 
$Y'_{sing}\tth .$

Consider $X$ as a smooth manifold, of real dimension $2\tth n.$ 
Then, for every $P\in X$ we have the almost-complex structure 
$J:T_{{}_P}X\to T_{{}_P}X$ on the real tangent space to $X$ at 
$P.$ Recall that $J$ is a real automorphism such that \tth $J^2=-id.$
Now, if $P\in\Delta$ then the tangent space $T_{{}_P}\Delta$ is a subspace 
of $T_{{}_P}X$ of real dimension $r\tth ,$ because $\sigma_h$ is a real-analytic 
embeding, and $T_{{}_P}\Delta +J(\tth T_{{}_P}\Delta\tth )$ is the 
{\sl smallest} complex subspace of $T_{{}_P}X$ which contains 
$T_{{}_P}\Delta\tth .$ Moreover, if this point $P$ of $\Delta$ is also 
non-singular for $Y',$ then
$$
T_{{}_P}\Delta +J(\tth T_{{}_P}\Delta\tth )\tth\subset\tth
T_{{}_P}Y'\tth\subset\tth T_{{}_P}X
$$
In fact, $Y'$ is a smooth manifold in a neighbrhood of $P,$ and 
\lq\lq\tth$T_{{}_P}Y'$\tth '' in the above relation means the real 
tangent space to $Y'$ at $P;$ in particular \tth $T_{{}_P}\Delta\subset
T_{{}_P}Y'.$ But $Y'$ is also, locally at $P,$ a complex submanifold of $X,$ hence 
$J(\tth T_{{}_P}Y'\tth )=T_{{}_P}Y'.$ Therefore, from \tth ${\rm dim}\tth Y'\leq n-p$ 
we get
\begin{equation}\label{notraverso}
{\rm dim}_{{}_\C}\bigl(\tth T_{{}_P}\Delta +J(\tth T_{{}_P}\Delta\tth )
\tth\bigr)\tth\leq\tth n-p
\end{equation}
for every $P\in\Delta$ which is non-singular for $Y'.$ If (\ref{notraverso})
is written in coordinates, it amounts to the vanishing of certain determinants,
whose entries are real-analytic functions of $P.$ Therefore, since 
(\ref{notraverso}) is satisfied at the points of a not empty open subset
of $\Delta\tth ,$ and $\Delta$ is connected, we conclude that (\ref{notraverso})
is satisfied at every point of $\Delta\tth .$

\medskip

Finally, we can forget the algebraic subset $Y',$ and {\sl define} the space 
$S^{\tth p,i}$ as consisting of the classes $\xi\in H^{\tth i}(X,\C )$ whose 
Poincar\'e dual $PD(\xi )$ can be represented by a singular cycle $\Gamma$
such that, for every singular simplex $\sigma$ in it, and for every 
$P\in\sigma (\Delta_r),$ the relation (\ref{notraverso}) is satisfied.
In the sequel we will say that such a $\Gamma$ is a {\sl good representative} 
for $PD(\xi )\tth .$

\medskip
 
As a first consequence of this definition and the above discussion we get
\begin{equation}\label{crono}
N^{\tth p}H^{\tth i}\tth X\tth\subset\tth
S^{\tth p,i}\tth\cap\tth H^{\tth i}\tth X
\end{equation}
The main result of the paper is 

\begin{thm}\label{finalm}\ {\sl For any integers $i\geq 0$ and $p\geq 0$ 
we have the equality}
\begin{equation}\label{artemide}
N^{\tth p}H^{\tth i}\tth X\tth =\tth
S^{\tth p,i}\tth\cap\tth H^{\tth i}\tth X
\end{equation} 
\end{thm} 

\medskip

Notice that condition (\ref{notraverso}) is given on the 
individual simplexes of $\Gamma .$ This will be the basic
ingredient to prove

\begin{prop}\label{rational}\ {\sl The spaces $S^{\tth p,i}$ are 
rational, namely there is a $\Q$-subspace $W\subset H^{\tth i}\tth X$ 
such that \ $S^{\tth p,i}\tth =\tth W\otimes_{{}_\Q}\C\tth .$} 
\end{prop}
\noindent 
By Galois descent, equality (\ref{artemide}) 
then implies \ $W=N^{\tth p}H^{\tth i}\tth X,$ and, a posteriori, the
meaning of the spaces $S^{\tth p,i}$ turns out to be
\begin{equation}\label{significato}
S^{\tth p,i}\tth =\tth N^{\tth p}H^{\tth i}\tth X\otimes_{{}_\Q}\C
\end{equation} 
In particular, this relation and (\ref{GHC0}) toghether imply
$$
S^{\tth p,i} \tth\subset\tth F^{\tth p}H^{\tth i}(X,\C )
$$
which shows that (\ref{crono}) is a refinement of (\ref{GHC0}). 

\medskip

These results leave untouched Grothendiek's Generalized Hodge Conjecture. 
On the other hand, it seems to me that the method described here to build 
algebraic supports for suitable cohomology classes is of independent 
interest. 

\medskip

After a section devoted to preliminary topics, the proofs of Theorem 
\ref{finalm} and Proposition \ref{rational} 
are given respectively in \S\S\ 3 and 4\tth . 

%%%%%%%%%%%%%%%%%%%%%%%%%%%%%%%%%%%%%%%%%%%%%%%%%%%%%%%%%%%%%%%%%%%%%%%%%%%
%%%%%%%%%%%%%%%%%%%%%%%%%%%%%%%%%%%%%%%%%%%%%%%%%%%%%%%%%%%%%%%%%%%%%%%%%%%

\section{\bf Preliminaries}
\label{premesse}  

Throughout this section we fix an embedding $X\subset\P^N$ and
a system of homogeneous coordinates
$(\tth\zeta_{\sss 0}:\zeta_{\sss 1}:\ldots :\zeta_{\sss N}\tth )$ 
for $\P^N.$

\medskip

First of all, we introduce a restriction on the systems of local 
holomorphic coordinates $(z_{\sss 1},\ldots ,z_n)$ on $X$ we will use 
in the sequel. 

Fix a point $P\in X;$ let 
$\Lambda\subset\P^N$ be any projective subspace of dimension $N-n-1\tth ,$ 
which is disjoint from $X$ and from the (\tth projective\tth ) tangent 
space to $X$ at $P.$ Finally, fix a general linear subspace $\P^n\subset\P^N.$
Let $\pi :X\to\P^n$ be the projection from $\Lambda\tth .$ 
Then the differential map $d\pi_{{}_P}$ of $\pi$ at $P$ 
is an isomorphism by construction. Hence the restriction of $\pi$ 
to a suitable euclidean neighborhood $V$ of $P$ is biholomorphic 
onto $\pi (V)\tth .$ 
{\sl We will use affine coordinates $z_{\sss 1},\ldots ,z_n$ on 
$\pi (V)$ as holomorphic coordinates on $V.$}

\smallskip

It is clear that these coordinates have an immediate algebro-geometricic meaning. 
In fact, they are essentially regular functions used cum grano salis. 
More precisely, assume that 
any point of a locus $L\subset V$ satisfies a system of algebraic 
equations $\varphi_j(z_{\sss 1},\ldots ,z_n)=0\tth ,$ with 
$j\in J\tth .$ These equations define an algebraic subset $Z$  
inside $\P^n,$ and $L$ is clearly contained into the intersection 
of $V$ with the cone $C$ projecting $Z$ from $\Lambda\tth .$ Moreover, 
if ${\rm codim}_{\pi (P)}(Z,\P^n)=t\tth ,$ then \ 
${\rm codim}_P(C\cap X,X)=t\tth .$ This follows easily from the fact that
$\pi$ is \'etale at $P.$

We will show at the end of \S\tth\ref{morehope} that the use of this kind of local
holomorphic coordinates is essential for the proof of Theorem \ref{finalm}.

\medskip

\noindent 
{\sl From now on, unless otherwise stated, 
we will use on $X$ only systems of local holomorphic coordinates 
as constructed above. Moreover we will make the harmless assumption that
the domain $V$ of every such system is contained in a Zariski open set  
$U_k\tth :=\tth\{\tth (\tth\zeta_{\sss 0}:\zeta_{\sss 1}:\ldots :
\zeta_{\sss N}\tth )\in\P^N\tth\vert\ \zeta_k\neq 0\tth\}$ of $\P^N,$
for some $k.$} 

\medskip
 
The rest of this section is devoted to the proof of the following lemma,
which in particular justifies the assumptions about the singular cycles to 
be used on $X$ we made in Introduction.

\begin{lemma}\label{omologia} Every class in $H_r(X,\Q )$ can be represented 
by a singular $r$-cycle $\Gamma = \sum_h\tth m_h\tth\sigma_h$ such that the
restriction to $\Delta_r$ of every singular simplex $\sigma_h:\bar\Delta_r\to X$ 
in it satisfies the following conditions
\begin{enumerate}
\item $\sigma_h:\Delta_r\to X$ is a real-analytic embedding;
\item the map $\sigma_h:\Delta_r\to X$ is semi-algebraic;
\item any $\Delta$ as above
is contained in a domain of local holomorphic coordinates on $X;$
hence, in particular, it is contained in some Zariski open set $U_k$ 
of $\P^N.$
\end{enumerate}
The $\Q$-vector space generated by such $r$-cycles will be denoted by
${\mathscr Z}_r\tth .$
\end{lemma}

\proof Fix a covering of $X$ by domains of holomorphic charts. Since the image 
of any closed singular simplex $\bar\Delta\to X$ is compact, by iterating the 
barycentric subdivision finitely many times we can replace such image by a chain,
all of whose simplexes are contained in some element of the covering. If $\Gamma$
is any singular $r$-cycle on $X,$ by applying this procedure to every simplex of 
$\Gamma$ we get a new $r$-cycle $\Gamma'$ on $X,$ homologous to $\Gamma ,$ which
satisfies condition 3. of the statement.

To get properties 1. and 2. we use the following theorem to construct a suitable 
triangulation of $X,$ and then we apply the isomorphism between simplicial 
and singular homology (see e.g. \cite{Mu}, \S\tth 34 and also p. 311). 
Before state it, let us recall the definition and a 
few, basic facts on semi-algebraic sets. The reader is referred 
to \cite{BCR} for a detailed account; a brief summary 
of their main properties is in \cite{Hi}.

Semi-algebraic sets are the subsets of some $\R^M$ which can be 
obtained by finite union, finite intersection and complementary set
starting from the family of sets 
$$
\{\tth x\in\R^M\tth\vert\tth f(x)\geq 0\tth\}
$$
where $f$ is a polinomial with real coefficients. Real algebraic 
subsets of $\R^M$ are elementary examples of this notion. 
The product of two semi-algebraic sets is still semi-algebraic.
A map between semi-algebraic sets is semi-algebraic if its graph is 
semi-algebraic. These maps are not necessarily continuous
(\tth however, we will deal only with continuous semi-algebraic maps\tth ).

\begin{thm}\label{trianghiro}\ {\rm \cite{Hi}}
{\sl Given a
finite system of bounded semi-algebraic sets $\{ X_\alpha\}$ in $\R^M,$ 
there exists a simplicial decomposition $\R^M=\cup_a\Delta_a$ and a
continuous semi-algebraic automorphism $\kappa$ of $\R^M$ such that
\begin{enumerate}
\item each $X_\alpha$ is a finite union of some of the $\kappa (\Delta_a 
)\tth ,$ namely the triangulation of $\R^M$ induces a triangulation of 
$X_\alpha$ for every $\alpha\tth ;$
\item $\kappa (\Delta_a)$ is a locally closed smooth real-analytic
submanifold of $\R^M$ and $\kappa$ induces a real-analytic isomorphism
$\Delta_a\simeq\kappa (\Delta_a)\tth ,$ for every $a\tth .$    
\end{enumerate}}
\end{thm}

\noindent
Notice that $\kappa$ is not a real-analytic map, only its restriction to any
(\tth open\tth ) simplex $\Delta_a$ of the simplicial decomposition of $\R^M$ is such.

\noindent
To apply this theorem to our set-up we need the following lemma, which
summarizes some material from \cite{vdW}.

\begin{lemma}\ {\sl There is a real-analytic embedding \ $\rho :\P^N_\C\to
\R^{(N+1)^2}$ \ which sends complex algebraic subsets of $\P^N$ to (\tth
necessarily compact\tth ) real algebraic subsets of $\R^{(N+1)^2}.$}
\end{lemma}

\proof For the reader convenience we give a sketch of proof.
Assume that the homogeneous coordinates
$(\tth\zeta_{\sss 0}:\zeta_{\sss 1}:\ldots :\zeta_{\sss N}\tth )$
for a point $P$ in $\P^N$ satisfy the normalization 
\begin{equation}\label{spaziopr}
\zeta_{\sss 0}{\bar\zeta}_{\sss 0}+\zeta_{\sss 1}{\bar\zeta}_{\sss 1}
+\ldots +\zeta_{\sss N}{\bar\zeta}_{\sss N}\tth =\tth 1
\end{equation}
Then the hermitian $(N+1)\times (N+1)$ matrix 
(\tth$\zeta_{\sss h}{\bar\zeta}_{\sss k}$\tth ) depends only on $P.$ 
Therefore the real and imaginary parts of the entries of this matrix
\begin{equation}\label{st}
\zeta_{\sss h}{\bar\zeta}_{\sss k}\tth =\tth
\sigma_{\sss hk}+i\tau_{\sss hk}\ \ \ \ \ \ \ \ 
0\leq h\leq k\leq N
\end{equation}
allow us to define the map $\rho$ by setting
$$
\rho :(\tth\zeta_{\sss 0}:\zeta_{\sss 1}:\ldots 
:\zeta_{\sss N}\tth )
\mapsto (\sigma_{\sss 00}\tth ,\sigma_{\sss 01},\ldots ,\sigma_{\sss NN},
\tau_{\sss 01},\ldots ,\tau_{\sss N-1,N})
$$
We leave to the reader to check that $\rho$ is a real-analytic embedding. 
The image $\mathfrak R$ of $\rho$ is compact, and homeomorphic to $\P^N.$ 
Moreover, $\mathfrak R$ is a real algebraic subset of $\R^{(N+1)^2}.$ 
In fact, a complete set of equations for $\mathfrak R$ is given by
$$
\sigma_{\sss 00}+\sigma_{\sss 11}+\ldots +\sigma_{\sss NN}\tth =\tth 1
$$
which translates (\ref{spaziopr}), and the equations obtained by separating
the real and imaginary part (\tth written in terms of the $\sigma$'s and 
$\tau$'s\tth ) from all the obvious relations\tth\footnote{
These relations amount to \ 
${\rm rk}\tth (\tth\zeta_{\sss h}{\bar\zeta}_{\sss k}\tth )\tth =1\tth .$}
\begin{equation}\label{licia}
\bigl(\tth\zeta_j{\bar\zeta}_k\tth\bigr)\bigl(\tth\zeta_u{\bar\zeta}_v
\tth\bigr)\tth -\tth\bigl(\tth
\zeta_j{\bar\zeta}_v\tth\bigr)\bigl(\tth\zeta_u{\bar\zeta}_k\tth\bigr)
\tth =\tth 0
\end{equation}

To complete the proof it remains to show that the image $\rho (Z)$ of any 
complex algebraic subset $Z\subset\P^N$ is a real algebraic subset of $\mathfrak 
R\tth .$ Let $\{ f_\nu\}_\nu$ be homogeneous generators for 
the saturated
homogeneous ideal of $Z\tth .$ Then \ $(\zeta )=(\tth\zeta_{\sss 0}:
\zeta_{\sss 1}:\ldots :\zeta_{\sss N}\tth )$ is in $Z$ iff $f_\nu (\zeta )=0$
for any $\nu\tth ,$ which is in turn equivalent to
$$
f_\nu (\zeta )\tth\bar f_\nu (\bar\zeta )\tth =\tth 0\ \ \ \ \ \ \ \ 
\ \ \ \ \ \ \ \ \hbox{for every}\ \ \nu
$$
Now, the left hand side of this relation can be written as a homogeneous polynomial 
with {\sl real} coefficients, in the indeterminates $\sigma_{\sss hk}$ and 
$\tau_{\sss hk}\tth .$ We will write it as $g_\nu\tth ,$ although it is
by no way uniquely determined. In fact, in any monomial 
$$
u\tth\zeta_{\sss 0}^{m_{\sss 0}}\tth\zeta_{\sss 1}^{m_{\sss 1}}\tth\ldots\tth
\zeta_{\sss N}^{m_{\sss N}}\tth\bar\zeta_{\sss 0}^{\mu_{\sss 0}}\tth\bar
\zeta_{\sss 1}^{\mu_{\sss 1}}\tth\ldots\tth\bar\zeta_{\sss N}^{\mu_{\sss N}}
\ \ \ \ \ \ \ \ \ \ \ \ \ \ m_{\sss 0}+\ldots +m_{\sss N}=
\mu_{\sss 0}+\ldots +\mu_{\sss N}=deg(f_\nu )
$$
of $f_\nu (\zeta )\tth\bar f_\nu (\bar\zeta )\tth ,$ a given $\zeta$ can be coupled 
with a $\bar\zeta$ in several differen ways. However, for a fixed $\nu\tth ,$ the 
various $g_\nu$ we get are all congruent modulo the ideal of  $\mathfrak R\tth .$

Then, any point of $\rho (Z)$ satisfies all the real algebraic equations \ 
$g_\nu=0\tth .$ Conversely, every point of $\mathfrak R$ which satisfies
all these equations is the image by $\rho$ of a point $(\zeta )\in\P^N$
such that $f_\nu (\zeta )=0$ for every $\nu\tth ,$ hence which is in
$Z\tth .$ Therefore, $\rho (Z)$ is a real algebraic subset of $\R^{(N+1)^2}$
and the proof is complete.
\qed

\medskip

\noindent
We have now the following situation, where a particular simplex of the 
triangulation of $\rho (X)$ has been highlighted.
\begin{equation}\label{triangolazione}
\objectmargin {0.6pc}
\diagram
\R^{(N+1)^2}\rto^{\sim}_\kappa & \R^{(N+1)^2} & \\
 & {\mathfrak R}\ar@{^{(}->}[u]\rto^{{\rho}^{-1}} & \P^N \\
\Delta_a\ar@{^{(}->}[uu]\rto_\kappa & 
\rho (X)\ar@{^{(}->}[u]\rto_{{\rho}^{-1}} 
& X\ar@{^{(}->}[u]
\enddiagram
\end{equation}
It is clear that all the singular simplexes $\sigma_a\tth :=\tth
{\rho}^{-1}\circ\kappa :\Delta_a\to X$ triangulate directly $X.$ 
These simplexes are real-analytic embeddings because of 

\begin{lemma}\ {\sl $\mathfrak R$ is a real-analytic manifold, of dimension
$2N,$ and \ $\rho^{-1} :\mathfrak R\to\P^N$ \ is a real-analytic map.}
\end{lemma}

\proof 
The first part of the statement is a straightforward consequence of the
previous lemma. In fact, we have seen that $\mathfrak R$ is real-algebraic
and smooth.

To prove that ${\rho}^{-1}$ is real-analytic, it is enough to check it locally.
Consider the Zariski open set $U_{\sss 0}=\{\tth\zeta_{\sss 0}\neq 0\tth\}$ of
$\P^N;$ then
$\rho (U_{\sss 0})$ is the open subset $\sigma_{\sss 00}\neq 0$ of 
$\mathfrak R\tth .$ We identify $U_{\sss 0}$ with $\C^N$ by the usual map   
$$
(\tth z_{\sss 1},\ldots ,z_{\sss N}\tth )\mapsto\Bigl(\tth
\frac{\scriptstyle 1}{a}:\frac{z_{\sss 1}}{a}:\ldots :\frac{z_{\sss N}}{a}
\tth\Bigr)\ \ \ \ \ \ \ \ \hbox{where}\ \ \ \ a^2\tth :=\tth 1+
z_{\sss 1}{\bar z}_{\sss 1}+\ldots +z_{\sss N}{\bar z}_{\sss N}
$$
Of course, the choice of homogeneous coordinates has been made to satisfy 
(\ref{spaziopr}). If $z_j=x_j+iy_j\tth ,$ when we compute $\rho (P)$ 
we get in particular for any $j=1,\ldots ,N$
$$
\sigma_{{\sss 0}j}+i\tth\tau_{{\sss 0}j}
\tth =\tth
\zeta_{\sss 0}\bar\zeta_j 
\tth =\tth
\frac{\bar z_j}{a^{\sss 2}}
\tth =\tth
\frac{x_j-i\tth y_j}{a^{\sss 2}}
\ \ \ \ \ \ \ \ \ \ \ \ \ \ \ \ \ 
 \hbox{and also}\ \ \ \ \ \ \ \ \ \ \ \ \ \ \ \ \ 
\sigma_{\sss 00}\tth =\tth\frac{\scriptstyle 1}{a^{\sss 2}}
$$
Since $x_j$ and $y_j$ are real, then the map $\rho^{-1}$ is given on 
$\rho (U_{\sss 0})$ by 
\begin{equation}\label{rho}
x_j\tth =\tth 
\frac{\sigma_{{\sss 0}j}}{\sigma_{\sss 00}}
\ \ \ \ \ \ \ \ \ \ \ \ \ \ \ \ \ \ \ \ 
y_j\tth =\tth 
-\tth\frac{\tau_{{\sss 0}j}}{\sigma_{\sss 00}}
\end{equation}  
Hence $\rho^{-1}$ is real-analytic on $\rho (U_{\sss 0})\tth .$ The same 
argument works for any open subset $\{\tth\zeta_k\neq 0\tth\}\subset\P^N.$
\qed

\medskip

To conclude the proof of Proposition \ref{omologia} we will deal with 
condition 2. For this, notice that 
if $Z\subset\P^N$ is closed, algebraic, then ${\mathfrak R}\backslash\rho (Z)$
is semi-algebraic. In particular, when $Z$ is the hyperplane $\zeta_k=0$
we will denote this set $W_k\tth ,$  namely $W_k=\rho (U_k)\tth .$ 

\medskip

Consider now a singular simplex \ $\sigma ={\rho}^{-1}\circ\kappa :
\Delta_r\to X$ as above, except that here the subscript \lq\tth r\tth ' 
simply denotes the dimension. 
To simplify notations we assume that the image $\Delta$ 
is contained in the Zariski open $U_0\tth .$
Hence we can factorize $\sigma$ in a slighty different way as
$$
\diagram
\Delta_r\rto^-\kappa & W_{\sss 0}\rto^-{{\rho}^{-1}} &
U_{\sss 0}=\R^{2N}\rto^-\pi & \R^{{\sss 2}n}
\enddiagram
$$
where $\pi$ is essentially a linear projection because of the
peculiar local holomorphic coordinates we use on $X,$ as explained in
the first part of this section.

\begin{lemma}\label{lemsemalg}\ {\sl The above map $\sigma :\Delta_r\to
\R^{{\sss 2}n}$ is semi-algebraic.}
\end{lemma}

\proof $\Delta_r$ is clearly a semi-algebraic set. The map $\kappa$ is 
semi-algebraic by Theorem \ref{trianghiro}, and linear projections 
like $\pi$ are easily seen to be semi-algebraic. Since the composition 
of semi-algebraic maps is still semi-algebraic \cite{BCR}, Prop. 2.2.6, 
it remains to show that \ $\rho^{-1}$ \ enjoys this 
property. 

We remarked above that $W_{\sss 0}$ is semi-algebraic. Moreover, 
by (\ref{rho}) the graph of $\rho^{-1}$ is contained inside the
algebraic subset $K$ of $W_{\sss 0}\times U_{\sss 0}=
W_{\sss 0}\times\C^N\tth ,$ defined by the equations
$$
\sigma_{\sss 00}\tth x_j-\sigma_{{\sss 0}j}\tth =\tth 0\ \ \ \ \ \ 
\ \ \ \sigma_{\sss 00}\tth y_j+\tau_{{\sss 0}j}\tth =\tth 0\ \ 
\ \ \ \ \ \ \ \ \ j=1,2,\ldots ,N
$$
Conversely, if $(u,v)\in K,$ then trivially $v=\rho^{-1}(u)$
because of the (\ref{rho}).
We conclude that $K$ is the graph of $\rho^{-1},$ and the proof is 
complete.
\qed 

Therefore, also the proof of Lemma \ref{omologia} is complete.

%%%%%%%%%%%%%%%%%%%%%%%%%%%%%%%%%%%%%%%%%%%%%%%%%%%%%%%%%%%%%%%%%%%%%%%%%%%%%
%%%%%%%%%%%%%%%%%%%%%%%%%%%%%%%%%%%%%%%%%%%%%%%%%%%%%%%%%%%%%%%%%%%%%%%%%%%%%

\section{\bf Proof of Theorem \ref{finalm}}
\label{morehope}

Let $\xi\in S^{\tth p,i}\tth\cap\tth H^{i}X\tth ,$ and 
let $\Gamma =\sum_h\tth m_h\tth\sigma_h\tth\in\tth{\mathscr Z}_r$ 
be a good representative for $PD(\xi )\tth .$  We want to construct 
an algebraic subset of $X,$ of codimension $\geq p\tth ,$ containing 
$\Gamma .$ For this, it is sufficient to construct for any singular 
simplex $\sigma_h$ of $\Gamma$ an algebraic subset of codimension $\geq p$
containing $\sigma_h(\Delta_r)\tth .$ 
 
\medskip

Let $\sigma:\Delta_r\to X$ be one of the simplexes of
$\Gamma ,$ and denote its image by $\Delta\tth .$ Recall that by Lemma
\ref{omologia} we have that 
$\Delta$ is contained in the domain $V$ of a system of local holomorphic 
coordinates $z_{\sss 1},\ldots ,z_n$ on $X.$ Therefore, if $u_{\sss 1},
u_{\sss 2},\ldots ,u_r$ are coordinates in $\R^r$ and $z_j=x_j+i\tth y_j$ 
as usual, the map $\sigma$ can be written in coordinates 
$$
x_j=f_j(u_{\sss 1},u_{\sss 2},\ldots ,u_r)\ \ \ \ \ \ \ 
y_j=f_{j+n}(u_{\sss 1},\ldots ,u_r)\ \ \ \ \ \ \ \ \ \ \ \ 
j=1,\ldots ,n
$$
where any \tth $f_k:\Delta_r\to\R$ \tth is real-analytic because $\sigma$
is. For future use we will denote the jacobian matrix of $\sigma$ by
\begin{equation}\label{amatriciana}
\left(\begin{array}{c}
{\bf A}\\
{\bf B}
\end{array}\right)
\end{equation}
where {\bf A} and {\bf B} are both $n\times r$ matrices, with real-analytic 
functions as entries. Note that the rank of (\ref{amatriciana}) is
$r$ because $\sigma$ is an embedding.

Then we define a real-analytic map \ $g:\Delta_r\to\C^n$ by setting
for any $j=1,\ldots ,n$
\begin{equation}\label{pongo}
z_j=g_j(u_{\sss 1},u_{\sss 2},\ldots ,u_r)\tth :=\tth
f_j(u_{\sss 1},\ldots ,u_r)+i\tth 
f_{j+n}(u_{\sss 1},\ldots ,u_r)
\end{equation}
Every $g_j$ is given locally at any fixed point $Q\in\Delta_r$ by a power
series in the $u$'s, with positive convergence radius. Hence, as it is well 
known, these power series allow to extend $g_j$ to a holomorphic function
$F_j\tth ,$ defined on a suitable neighborhood of $\Delta_r$ inside $\C^r$ 
(\tth on $\C^r$ we will 
use complex coordinates $w_{\sss 1}, w_{\sss 2},\ldots ,w_r$ where 
$w_k=u_k+iv_k$ is the decomposition of $w_k$ into its real and 
imaginary parts, and the $u_k$'s are the old real coordinates of 
$\R^r$\tth ). Therefore, there is a {\sl connected}
open set $E\subset\C^r$ containing $\Delta_r\tth ,$
such that every $F_j$ is defined on $E,$ and the $F_j$
are the components of a holomorphic map $F:E\to\C^n.$ 
Moreover, by continuity of $F$ we can also assume (\tth with a slight 
abuse of notation\tth ) that $F(E)\subset V.$ 

A function like $F_j$ is usually called a {\sl holomorphic extension} 
of $g_j\tth .$ Its germ along $\Delta_r$ is unique. We will use this 
procedure later, preferably denoting by $\widetilde\varphi$ a 
holomorphic extension of a real-analytic function $\varphi\tth .$ 
This explains the notation used in the following statement.

\begin{lemma}\label{rangorozzo}\ {\sl The jacobian matrix $J$ of $F$ is
$\widetilde{\bf A}+i\tth\widetilde{\bf B}$ (\tth see (\ref{amatriciana})\tth ),
its rank is \ $\leq n-p$ \ at any point of $E,$ and reaches 
its maximum at some point of $\Delta_r\tth .$}
\end{lemma}

\proof
The first part is a straightforward consequence of the fact that 
clearly $J={\bf A}+i\tth{\bf B}$ along $\Delta_r\tth .$

Let $Z\subset E$ denote the set where the rank of $J$ is not maximum.
Then $Z$ is a complex analytic set, closed in $E,$ of dimension $<r\tth .$
Assume $\Delta_r\subset Z\tth .$ 

Let $\varphi :W\to\C$ be any holomorphic
function that vanishes on $W\cap Z\tth ,$ where $W\subset E$ is open, 
connected and $W\cap\Delta_r\neq\emptyset\tth .$ Then $\varphi$ 
vanishes identically on $W\cap\Delta_r\tth ,$ hence it vanishes identically 
on $W$ because $W$ is connected and \cite{S}, pag. 21, a variant of the Principle 
of Analytic Continuation. This implies that $Z$ has no non-trivial equations,
and therefore the dimension of $Z$ is $r\tth ,$ contradiction. 
Then $\Delta_r\not\subset Z$ and the proof is complete.
\qed

\medskip

Another basic property of the map $F:E\to\C^n$ \ is

\begin{lemma}\label{finitezza}\ 
{\sl Let $E_{\sss 0}\subseteq E$ be any open, semi-algebraic subset. Then the 
restriction of $F$ to $E_{\sss 0}$ is a semi-algebraic map
$F:E_{\sss 0}\to\C^n=\R^{{\sss 2}n}.$ 
} 
\end{lemma}

\proof
We know by Lemma \ref{lemsemalg} that $\sigma :\Delta_r\to\R^{{\sss 2}n}$ 
is semi-algebraic. It is easily seen that this is equivalent to have
$f_j:\Delta_r\to\R$ semi-algebraic for every $j=1,\ldots ,2n\tth .$
By \cite{BCR}, Prop. 8.1.8, since any $f_j$ is of class ${\mathscr C}^\infty,$ 
this amounts to the existence for any $j$ of a polynomial 
$P_j(U_{\sss 1},\ldots ,U_r,T)$ with real coefficients and positive 
degree with respect to $T,$ such that for every $u\in\Delta_r$ 
\begin{equation}\label{basta}
P_j(u,f_j(u))\tth =\tth 0
\end{equation}
holds true. The left hand side 
of the above relation can be thought as a real-analytic
function of $u\tth .$ Since $E$ has been assumed to be connected, 
if we replace $u$ by $u+iv$ in 
(\ref{basta}) we get by \cite{S}, pag. 21, that
\begin{equation}\label{hasta}
P_j(u+iv,f_j(u+iv))\tth =\tth 0
\end{equation}
for every $u+iv\in E$.

This relation says that $f_j(u+iv)$ is algebraic over the field 
of fractions $\C (u,v)$ of the ring of polynomials $\C [u,v]=
\C [u_{\sss 1},\ldots ,u_r,v_{\sss 1},\ldots ,v_r]\tth .$
But $\C (u,v)$ is algebraic over $\R (u,v)\tth ,$ hence $f_j(u+iv)$ 
is algebraic over $\R (u,v)$ as well, and there is 
$Q_j\in\R [u,v,T]\tth ,$ with positive degree with respect to 
$T,$ such that for every $(u,v)=u+iv\in E$ 
$$
Q_j(u,v,f_j(u+iv))\tth =\tth  0
$$
Since $Q_j$ has real coefficients, this implies \ $Q_j(u,v,
\overline{f_j(u+iv)}\tth )=0\tth ,$ hence $\overline{f_j(u+iv)}$ 
is also algebraic over $\R (u,v)\tth .$ If we set 
$$
f_j(u+iv)=\varphi_j(u+iv)+i\tth\psi_j(u+iv)
$$ 
where $\varphi_j$ and $\psi_j$ are real valued functions
defined on $E,$ then we can conclude that both $\varphi_j(u+iv)$
and $\psi_j(u+iv)$ are algebraic over $\R (u,v)\tth .$ 
Therefore these functions are
semi-algebraic on every open, semi-algebraic subset $E_{\sss 0}$ of 
$E\subset R^{\tth{\sss 2}r}$ \cite{BCR}, Prop. 8.1.8.
Finally, by (\ref{pongo}) we have
$$
\begin{array}{l}
F_j(u+iv)\tth =\tth f_j(u+iv)+i\tth f_{j+n}(u+iv)\tth =\tth \\
\\
=\tth \varphi_j(u+iv)-\psi_{j+n}(u+iv)+i\tth\bigl(
\psi_j(u+iv)+\varphi_{j+n}(u+iv)\bigr)
\end{array}
$$
Hence, for any $j=1,\ldots ,n\tth ,$ both the real and imaginary
components of $F_j$ are semi-algebraic functions on $E_{\sss 0}\tth ,$ 
and the proof of Lemma \ref{finitezza} is complete.
\qed

\bigskip

Let $P_{\sss 0}\in\Delta_r$ be a point where the rank of the Jacobian matrix 
of $F$ reaches its maximum. Then by the holomorphic variant of the Rank Theorem 
(\tth see e.g. \cite{KK}\tth ),
the map $F$ can be written locally at $P_{\sss 0}$ as a linear 
projection, with respect to suitable local holomorphic coordinates 
(\tth not necessarily of the kind introduced in \S\tth\ref{premesse}\tth ). 
This implies that $F$ is an open map, locally at $P_{\sss 0}\tth ,$ and that a
neighborhood of $P=F(P_{\sss 0})$ on $F(E)$ is a complex-analytic manifold. 
More precisely,
there is a suitably small open polydisk $W\subset E,$ centered at 
$P_{\sss 0}\tth ,$ such that $Y:=F(W)$ is a complex manifold, locally closed 
in $V.$ The dimension of $Y$ is $n-t\leq n-p\tth ,$ by Lemma \ref{rangorozzo}.
Moreover, $Y$ is semi-algebraic by Lemma \ref{finitezza} and \cite{BCR}, 
Prop. 2.2.7.

\medskip

To prove Theorem \ref{finalm} it is then sufficient to check that the
dimension of the Zariski closure $Z$ of $Y$ inside $X,$ with respect to the 
{\sl complex} Zariski topology, is equal to the dimension of the complex 
manifold $Y,$ and therefore is $\leq n-p\tth .$ In fact, by construction of 
$Y,$ a non empty open subset of $\Delta_r$ is sent by $\sigma$ into $Z,$ 
hence $\sigma (\Delta_r)\subset Z$ because $\sigma$ is real-analytic.

In the beautyful paper \cite{C} it was recognized that the property for 
a complex-analytic subset $Y$ of some algebraic complex manifold to be itself 
algebraic is of {\sl local} nature. Heuristically, this suggests that to show
${\rm dim}\tth Z\leq n-p$ we can argue in the ring $\cal O$ of germs of
holomorphic functions on $X$ at a suitable point $Q$ of $Y.$ 

For the proof, a last ingredient is needed, namely that $Y$ is a semi-algebraic set.
That semi-algebraicity can be useful to deal with our problem is showed by 
the fact that a complex-analytic subset $Y$ of $\C^m$ which is semi-algebraic 
as subset of $\R^{{\sss 2}m}=\C^m$ is actually complex algebraic \cite{FL}.
After the paper was completed, I learned that local versions of this
result were obtained in \cite{FLR}.

\medskip

Denote by $\bar Y$ the Zariski closure of $Y$ inside $\R^{{\sss 2}n}$
(\tth recall that $F:E\to\C^n=\R^{{\sss 2}n}$\tth ). The key argument in the 
following pages is a local computation at a point of $Y.$ To perform it, we
need that this point is non singular for $\bar Y.$ I was unable to prove
(\tth or disprove\tth ) this for the point $P=F(P_{\sss 0})$ introduced above. 
However

\begin{lemma}\label{fabrizio}\ 
{\sl The algebraic set $\bar Y$ is irreducible and $Y$
contains points which are non singular for $\bar Y.$}
\end{lemma}

\proof
Let $A$ denote the ring of germs at $P$ of real-analytic, complex valued 
functions on $X.$ Moreover, denote by ${\mathscr J}\subset A$ the ideal of 
germs of functions whose restriction to $Y$ vanishes identically around $P.$ 
Since $Y$ is also a real-analytic manifold, of dimension $2(n-t)\tth ,$ at 
every
point it is locally ${\mathscr C}^\omega$-diffeomorphic to some open subset
of $\R^{{\sss 2}(n-t)}.$ Hence the germ of $Y$ at $P$ is irreducible. A
standard argument then shows that $\mathscr J$ is a prime ideal.

Since the local coordinates at $P$ are fixed, there is a subring of $A$
which is isomorphic to the polynomial ring 
$\R [\tth\underline x\tth ,\underline y\tth ]\tth .$  
But the ideal of $\bar Y$ is clearly ${\mathscr J}\cap\R [\tth\underline 
x\tth ,\underline y\tth ]\tth ,$ hence $\bar Y$ is irreducible.

\medskip

Now, denote by $\bar Y_{reg}$ and $\bar Y_{sing}$ respectively the sets of 
non-singular and singular points of $\bar Y.$ If $Y\cap\bar Y_{reg}=\emptyset\tth ,$ 
then $Y\subset\bar Y_{sing}\subsetneq\bar Y.$ Since $\bar Y_{sing}$ is 
a closed algebraic subset of $\R^{{\sss 2}(n-t)},$ we have a contradiction.
\qed

\medskip

From now on $Q$ will denote a fixed point of $Y,$ which is non singular for 
$\bar Y.$ To simplify notations, we will also assume that local holomorphic 
coordinates have been fixed around $Q$ so that $Q=(0,\ldots ,0)\tth .$ 

\medskip

Let $\cal O$ denote the ring of germs at $Q$ of holomorphic functions
on $X,$ and let ${\mathscr I}\subset{\mathcal O}$ be the ideal of germs 
representing functions whose restriction to $Y$ vanishes identically 
around $Q.$ We will see in a moment that $\mathscr I$ is a prime ideal.
Moreover, ${\rm ht}({\mathscr I})=t$ because the dimension of $Y$ as
complex manifold is 
$n-t\tth .$ As in the proof of Lemma \ref{fabrizio}, 
any choice of a system of local holomorphic coordinates for $X$ at $Q$ determines 
a subring of $\cal O\tth ,$ which is isomorphic to the ring of polynomials 
$\C [z_{\sss 1},\ldots ,z_n]\tth ;$ hence this subring is not intrinsic. However, 
since we use on $X$ only the systems of local 
holomorphic coordinates introduced in \S\tth\ref{premesse}\tth ,
the elements of the corresponding subring 
$\C [z_{\sss 1},\ldots ,z_n]$ of $\cal O$ 
actually represent germs at $Q$ of {\sl regular} functions on $X.$ 
If we consider the (\tth prime\tth ) ideal
$I:={\mathscr I}\cap\C [\tth\underline z\tth ]\tth ,$ this implies that 
to conclude the proof of Theorem \ref{finalm} we have only to check that
\begin{equation}\label{numgiust}
{\rm ht}(I)\tth\geq\tth t
\end{equation}
(\tth it is not hard to show that also ${\rm ht}(I)\leq t$ holds true, 
so that we have an equality\tth ). In fact, by the above remarks, 
the Zariski closure $Z$ of
$Y$ for the complex Zariski topology of $X$ is defined locally at $Q$ by 
the ideal $I,$ and (\ref{numgiust}) implies that
$Z$ has codimension $\geq t\geq p$ in $X.$  

\medskip

To use the additional information that $Y$ is semi-algebraic, we have to 
complete the algebraic set-up introduced above. As in the proof of Lemma 
\ref{fabrizio}, let $A$ denote the ring of germs at $Q$ of real-analytic, complex 
valued functions on $X.$ Notice that $\mathcal O$ is a subring of 
$A\tth .$ Both $A$ and $\mathcal O$ are regular local rings, of dimension 
$2\tth n$ and $n$ respectively. Moreover,
$\C [\tth\underline x\tth ,\underline y\tth ]=\C 
[x_{\sss 1},\ldots ,x_n,y_{\sss 1},\ldots ,y_n]$ is a subring of $A\tth .$
Hence, to start we have the commutative diagram of rings
\begin{equation}\label{rings}
\diagram
{\mathcal O}\rto & A \\
\C [\tth\underline z\tth ]\uto\rto & \C [\tth\underline x\tth ,
\underline y\tth ]\uto
\enddiagram
\end{equation}
where all the maps are inclusions. Now, denote by ${\mathscr J}\subset A$ 
the ideal of germs of functions whose restriction to $Y$ vanishes identically 
around $Q\tth ,$ and set $J:={\mathscr J}\cap\C [\tth\underline x\tth ,
\underline y\tth ]\tth .$ The proof of Lemma \ref{fabrizio} shows that 
$\mathscr J$ is a prime ideal. Since ${\mathscr I}={\mathscr J}\cap
{\mathcal O}\tth ,$ \ $\mathscr I$ is a prime ideal as well.

\bigskip

To prove (\ref{numgiust}), we start by applying the \lq\tth dimension
formula\tth ' \cite{M}, p.\tth 119 to the
ring extension $\C [\tth\underline z\tth ]\subset \C [\tth\underline x
\tth ,\underline y\tth ]\tth ,$ thus getting (\tth notice that 
$\C [\tth\underline x\tth ,\underline y\tth ]=\C [\tth\underline z\tth ,
\underline x\tth ]$\tth )
\begin{equation}\label{dimfor1}
{\rm ht}(J)\tth +\tth{\rm tr.deg.}\tth_{{}_{\kappa(I)}}\tth\kappa (J)
\tth =\tth {\rm ht}(I)\tth +\tth 
{\rm tr.deg.}\tth_{{}_{\C [\tth\underline z\tth ]}}\tth 
\C [\tth\underline z\tth ,\underline x\tth ] 
\end{equation}
where the trascendence degree in the right hand side is that of the quotient 
field of $\C [\tth\underline x\tth ,\underline y\tth ]$ over that of
$\C [\tth\underline z\tth ]\tth ,$ and $\kappa (J)$ is the quotient field
of $\C [\tth\underline x\tth ,\underline y\tth ]/J\tth .$ 
\begin{lemma}\label{bu} 
${\rm ht}(J)\tth =\tth 2\tth t$ 
\end{lemma}

\proof The real-analytic manifold underlying $Y$ has dimension 
$2(n-t)\tth ,$ and this is also the dimension of $Y$ as a semi-algebraic 
set \cite{BCR}, Prop. 2.8.14.

The ideal $K\subset\R [\tth\underline x\tth ,
\underline y\tth ]$ of all the polynomials vanishing on $Y$ 
defines the Zariski closure $\overline Y$ of $Y$ inside 
$\R^{\tth{\sss 2}n}.$ It is clear that $K=J\cap\R [\tth
\underline x\tth ,\underline y\tth ]\tth ,$ which implies that $K$ 
is a prime ideal. By \cite{BCR}, Prop. 2.8.2, the dimension 
of \tth $\overline Y$ \tth is $2(n-t)\tth ,$ and then \tth 
${\rm ht}(K)=2\tth t$ \cite{M2}, (14.H),p. 92.

Finally, since the extension $\R [\tth\underline x\tth ,\underline 
y\tth ]\subset\C [\tth\underline x\tth ,\underline y\tth ]$ is
integral and flat, we can conclude ${\rm ht}(J)={\rm ht}(K)=2\tth t$
\ \cite{M2}, (13.C), p. 81.
\qed

\medskip

Replacing ${\rm ht}(J)=2\tth t$ into (\ref{dimfor1}) shows 
that (\ref{numgiust}) is equivalent to
\begin{equation}\label{numgiust2}
{\rm tr.deg.}\tth_{{}_{\kappa(I)}}\tth\kappa (J)\tth \geq \tth n-t
\end{equation}

\medskip

We will apply now the \lq\tth dimension formula\tth ' again, this
time to the ring extension ${\mathcal O}\to{\mathcal O}[\tth
\underline x\tth ]$ (\tth $\subset A$\tth ). This is possible because 
${\mathcal O}$ is universally catenarian, being a regular local ring, 
and ${\mathcal O}[\tth\underline x\tth ]$ is an 
$\mathcal O$-algebra of finite type. Set
${\mathscr H}:={\mathscr J}\cap{\mathcal O}[\tth
\underline x\tth ]\tth .$ Then 
$$
{\rm ht}({\mathscr H})\tth +\tth 
{\rm tr.deg.}\tth_{{}_{\kappa ({\mathscr I})}}\tth\kappa ({\mathscr H})
\tth =\tth
{\rm ht}({\mathscr I})\tth +\tth 
{\rm tr.deg.}\tth_{{}_{{\mathcal O}}}{\mathcal O}[\tth\underline x\tth ]
$$
The rings $\mathcal O$ and $A$ can be thought as the rings of convergent 
power series with complex coefficients, respectively in the variables
$z_{\sss 1},\ldots ,z_n$ and $z_{\sss 1},\ldots ,z_n,x_{\sss 1},\ldots 
,x_n\tth .$ This shows that $x_{\sss 1},\ldots 
,x_n$ are algebraically independent over $\mathcal O\tth .$ 
Hence, if we assume for a moment 
\begin{lemma}\label{bah} 
${\rm ht}({\mathscr H})\tth =\tth 2\tth t$
\end{lemma}

\noindent
we conclude \ ${\rm tr.deg.}\tth_{{}_{\kappa({\mathscr I})}}\tth
\kappa ({\mathscr H})= n-t\tth .$

\medskip

Finally, consider the commutative diagram of integral domains and 
injective ring maps, and the corresponding diagram of 
quotient fields
$$
\diagram
{\mathcal O}/{\mathscr I}\rto & {\mathcal O}[\tth\underline x\tth ]/
{\mathscr H} \\
\C [\tth\underline z\tth ]/I\uto\rto & \C [\tth\underline z\tth ,
\underline x\tth ]/J\uto
\enddiagram
\ \ \ \ \ \ \ \ \ \ \ \ \ 
\diagram
\kappa ({\mathscr I})\rto & \kappa ({\mathscr H}) \\
\kappa (I)\uto\rto & \kappa (J)\uto
\enddiagram
$$
It is clear that a trascendence base for $\kappa ({\mathscr H})$ over
$\kappa ({\mathscr I})$ can be extracted from the set of the residue 
classes mod $\mathscr H$ of $x_{\sss 1},\ldots ,x_n\tth .$ Say 
such a base is $\bar x_{\sss 1},\ldots ,\bar x_{n-t}\tth ,$ the 
number of its elements was determined above. It is also clear that 
these elements are, a fortiori, algebraically independent over 
$\kappa (I)\tth .$ But $\bar x_{\sss 1},\ldots ,\bar x_{n-t}$ actually
belong to $\C [\tth\underline z\tth ,\underline x\tth ]/J,$ and
the proof of (\ref{numgiust2}) is complete, except for the proof of
Lemma \ref{bah}. 

\bigskip

\dimost of Lemma \ref{bah}. First of all, we construct 
a suitable system of generators for the ideal $J_{loc}\subset
\C [\tth\underline z\tth ,\underline x\tth ]_{loc}\tth ,$
where \lq\tth$loc$\tth ' denotes the localization
with respect to $(\tth\underline z\tth ,\underline x\tth )\tth .$
Recall that we assumed $Q=(0,\ldots ,0)\tth ;$ 
we can also assume that $Y$ is given locally at $Q$ by  
holomorphic equations like
$$
z_{\sss 1}\tth +\tth\hbox{higher order terms}\tth =\tth 0
\ \ \ \ \ \ 
\ldots\ldots \ \ \ z_t\tth +\tth\hbox{h.o.t.}\tth =\tth 0
$$
Hence, if we consider $Y\subset V\subseteq\C^n,$ then the 
embedded tangent space to the complex manifold $Y$ at $Q$ 
is defined by \ $z_{\sss 1}=\ldots =z_t=0\tth .$ Therefore 
equations for the embedded tangent space to the underlying 
real-analytic manifold $Y\subset\R^{{\sss 2}n}$ are
$$
x_{\sss 1}\tth =\tth 0\ \ \ \ldots\ldots\ \ \ x_t\tth =\tth 0\ \ 
\ \ y_{\sss 1}\tth =\tth 0\ \ \ \ldots\ldots\ \ \ y_t\tth =\tth 0
$$
As it is remarked in \cite{BCR}, p. 66, after Definition 3.3.4,
since $Q$ is a non singular point of $\bar Y,$ then the above equations 
define also the Zariski tangent space to the algebraic set $\bar Y$
at $Q\tth .$
Then, with the notations of the proof of Lemma \ref{bu}, there are 
elements $P_{\sss 1},\ldots , P_{{\sss 2}t}\in K\subset\R [\tth
\underline x\tth ,\underline y\tth ]$ such that for any $j=1,\ldots ,t$
$$
P_j\tth =\tth x_j\tth +\tth\hbox{h.o.t.}
\ \ \ \ \ \ldots\ldots\ \ \ \ \  
P_{t+j}\tth =\tth y_j\tth +\tth\hbox{h.o.t.}
$$
and the equations
\begin{equation}\label{equa2}
P_j\tth =\tth 0
\ \ \ \ \ \ldots\ldots\ \ \ \ \ 
P_{t+j}\tth =\tth 0
\end{equation}
define the semi-algebraic set $Y$
in a suitable euclidean neighborhood of $Q.$ 

These polynomials $P_h$ are contained in the maximal ideals of the local
rings $\C [\tth\underline z\tth ,\underline x\tth ]_{loc}$ and $A\tth .$
More precisely, these rings are regular local rings and $P_{\sss 1},\ldots , 
P_{{\sss 2}t}$ is a subset of a minimal system of parameters for both. Therefore, 
by \cite{M}, theorems 14.2 and 14.3, we get the following conclusion,
that we record as a lemma for easy of reference.

\begin{lemma}\label{parametri}\ {\sl $P_{\sss 1},\ldots , 
P_t,P_{t+{\sss 1}}\tth ,\ldots ,P_{{\sss 2}t}$ 
generate a prime ideal 
of height $2\tth t$ in both rings $\C [\tth\underline z\tth ,\underline 
x\tth ]_{loc}$ \ and \ $A\tth .$
}
\end{lemma}

\noindent
Now the ideal $(P_{\sss 1},\ldots ,
P_{{\sss 2}\tth t})$ is clearly contained into $J_{loc}\tth ,$ hence the two
ideals coincide having the same height, and we have the desired system of 
generators for $J_{loc}\tth .$ 

\medskip

\noindent
{\bf Claim} The polynomials $P_{\sss 1},\ldots ,P_{{\sss 2}t}$ 
generate the ideal $\mathscr J$ in $A\tth .$

\medskip

\noindent
Clearly $\mathscr J$ contains these polynomials.   Conversely,
let $\varphi :T\to\C$ a real-analytic function, with $T\subset
\R^{{\sss 2}n}$ a neighborhood of $Q$ in $X,$ such that its germ 
belongs to $\mathscr J.$ By definition of $\mathscr J$ this
means that the restriction of $\varphi$ to $Y$ vanishes identically
around $Q.$

Consider now the algebraic subset $D$ of $\C^{{\sss 2}n},$ defined by the 
equations (\ref{equa2}). Locally at $\underline 0 =(0,\ldots ,0)\in
\C^{{\sss 2}n}$ it is a complex manifold, of dimension $2\tth n-2\tth t
\tth .$
We already remarked (\tth just before Lemma \ref{bah}\tth ) that
$A$ can be thought as the ring of convergent power series with complex 
coefficients, in the variables $z_{\sss 1},\ldots ,z_n,x_{\sss 1},
\ldots ,x_n\tth .$ Hence $A$ can be identified with the ring of germs at
$\underline 0$ of holomorphic functions $\C^{{\sss 2}n}\to\C\tth .$
The holomorphic extension $\widetilde\varphi$ of $\varphi$ is the 
holomorphic function defined in a suitable neighborhood of $\underline 0$
inside $\C^{{\sss 2}n}$ by the same power series of $\varphi\tth .$
Therefore, to prove that the germ of $\varphi$ belongs to the ideal of $A$
generated by the polynomials $P_{\sss 1},\ldots ,P_{{\sss 2}t}\tth ,$ 
we have only to check that the restriction 
of \tth $\widetilde\varphi$ \tth to $D$ vanishes identically around 
the point $Q$
because of the analytic version of Hilbert'\tth s Nullstellensatz \cite{GR}, 
Ch. II, Section E (\tth in particular
Thm. 11 p. 88; see also p. 90\tth ). In fact, by Lemma \ref{parametri}
the polynomials $P_{\sss 1},\ldots ,P_{{\sss 2}t}$ generate a prime ideal 
in the ring $A\tth ,$ hence
$$
\sqrt{(P_{\sss 1},\ldots ,P_{{\sss 2}t})\tth A}
\tth =\tth (P_{\sss 1},\ldots ,P_{{\sss 2}t})\tth A
$$
If we apply the Implicit Function Theorem for real-analytic functions to 
(\ref{equa2}), we get a neighborhood $H$ of $\underline 0\in\R^{{\sss 2}n-
{\sss 2}t},$ with coordinates $x_{t+{\sss 1}}\tth ,\ldots ,x_n,
y_{t+{\sss 1}}\tth ,\ldots ,y_n\tth ,$ 
and a real-analytic map $\alpha :H\to\R^{{\sss 2}n}$ 
which parametrizes $Y$ locally at $Q.$ Therefore the holomorphic
extension $\widetilde\alpha$ of $\alpha$ parametrizes $D$ locally at 
$\underline 0\tth .$

The holomorphic function $\widetilde\varphi\circ\widetilde\alpha -
\widetilde{\varphi\circ\alpha}$ is defined in a neighborhood $H'$ of $H$
inside $\C^{{\sss 2}n-{\sss 2}t};$ if $H$ is connected, say an open ball,
it is harmless to assume that $H'$ is connected as well. But the
restriction of $\widetilde\varphi\circ\widetilde\alpha -
\widetilde{\varphi\circ\alpha}$ to $H$ vanishes identically, hence
$\widetilde\varphi\circ\widetilde\alpha =\widetilde{\varphi\circ\alpha}$
on $H'$ \cite{S}, pag. 21. Now, since $\varphi\equiv 0$ on $Y,$ we have 
$\varphi\circ\alpha\equiv 0\tth ,$ hence \tth $\widetilde\varphi
\circ\widetilde\alpha\equiv 0\tth .$ In other words,
the restriction of \tth $\widetilde\varphi$ \tth to $D$ 
vanishes identically around $\underline 0\tth .$ As we said above, this 
implies that the germ of $\varphi$ belongs to $(P_{\sss 1},\ldots ,
P_{{\sss 2}t})\tth A\tth ,$ and therefore the Claim is completely proved.

\medskip

To complete the proof of Lemma \ref{bah},
let $\mathscr M$ denote the maximal ideal of $\mathcal O\tth .$
Then ${\mathscr N}:={\mathscr M}+(\tth\underline x\tth )$ is a maximal
ideal of ${\mathcal O}[\tth\underline x\tth ]$ containing
${\mathscr H} ,$ and we have a local 
homomorphism ${\mathcal O}[\tth\underline x\tth ]_{\mathscr N}\to
A\tth .$ Both rings are regular local rings of dimension $2\tth n,$ 
whose respective maximal ideals are both generated by the germs of 
the functions $z_{\sss 1},\ldots ,z_n, x_{\sss 1},\ldots 
,x_n\tth .$ We then conclude that ${\mathcal O}[\tth\underline x\tth 
]_{\mathscr N}\to A$ is flat \cite{M}, Thm. 23.1, hence it is a faithfully 
flat extension because the homomorphism ${\mathcal O}[\tth\underline x\tth 
]_{\mathscr N}\to A$ is local and \cite{M}, Thm. 7.2.
But then, by the Claim and \cite{M}, Thm. 7.5
$$
(P_{\sss 1},\ldots ,P_{{\sss 2}t})\tth{\mathcal O}[\tth\underline 
x\tth ]_{\mathscr N}\tth =\tth 
(P_{\sss 1},\ldots ,P_{{\sss 2}t})\tth A\tth\cap\tth
{\mathcal O}[\tth\underline x\tth ]_{\mathscr N}\tth =\tth 
{\mathscr J}\tth\cap\tth
{\mathcal O}[\tth\underline x\tth ]_{\mathscr N}\tth =\tth
{\mathscr H}_{\mathscr N}
$$
The same argument used to prove Lemma \ref{parametri}
shows that $(P_{\sss 1},\ldots ,P_{{\sss 2}t})\tth
{\mathcal O}[\tth\underline x\tth ]_{\mathscr N}$ is a prime ideal of
height $2\tth t,$ and Lemma \ref{bah} is completely proved, 
as well as Theorem \ref{finalm}.

\bigskip

\noindent
{\bf Remark}
The use of the kind of holomorphic coordinates introduced in \S\tth\ref{premesse}
is {\sl necessary} for the above argument to work, as the following example shows.

Let $X$ be an abelian variety, of dimension $g\geq 2\tth ,$ and let $\pi :\C^g\to X$ 
be the universal covering of $X.$ The map $\pi$ is holomorphic and locally 
biholomorphic on the image. Then we can use cartesian complex coordinates 
$z_{\sss 1},\ldots ,z_g$ on $\C^g$ as local holomorphic coordinates on $X.$ 
Now, let $v\in\C^g$ be a non zero vector such that $\pi (\{\tth tv\vert t\in\R\tth\} )$ 
is dense in $X$ for the euclidean topology. The map
$$
\sigma :\R\to X\ \ \ \ \ \ \ \ \hbox{given by}\ \ \ \ \ \ \ \ 
\sigma (t)\tth :=\tth\pi (tv) 
$$
is real-analytic. Moreover, the restriction of $\sigma$ to $(0,1)$ is a real-analytic 
embedding, and if $\parallel\!\! v\!\parallel$ is sufficiently small, it is also immediately 
checked by using the local coordinates $z_{\sss 1},\ldots ,z_g$ that such restriction of
$\sigma$ is semi-algebraic.
Nevertheless, $\sigma ((0,1))$ {\sl is not contained in any complex algebraic subvariety}
$Z\subsetneq X.$ In fact, otherwise the real-analyticity of the map $\sigma :\R\to X$
would imply by a standard argument that $\sigma (\R )\subset Z,$ hence \ $\overline{\pi 
(\{\tth tv\vert t\in\R\tth\} )}\subsetneq X,$ contradiction by density.

This different behaviour of the local holomorphic coordinates $z_{\sss 1},\ldots ,z_g$
is easily understood if we think at $\pi$ as a map $\pi :\C^g\to\P^N,$ defined by suitable
theta functions $\theta_{\sss 0},\theta_{\sss 1},\ldots ,\theta_N\tth ,$ with image $X$ 
$$
\pi (z_{\sss 1},\ldots ,z_g)\tth =\tth
(\tth\theta_{\sss 0}(z_{\sss 1},\ldots ,z_g):\theta_{\sss 1}(z_{\sss 1},\ldots ,z_g):
\ldots :\theta_N(z_{\sss 1},\ldots ,z_g)\tth )
$$
Now the theta functions are trascendental, whereas the holomorphic coordinates 
introduced in \S\tth\ref{premesse} are essentially algebraic.

%%%%%%%%%%%%%%%%%%%%%%%%%%%%%%%%%%%%%%%%%%%%%%%%%%%%%%%%%%%%%%%%%%%%%%%%%%%%%%
%%%%%%%%%%%%%%%%%%%%%%%%%%%%%%%%%%%%%%%%%%%%%%%%%%%%%%%%%%%%%%%%%%%%%%%%%%%%%%

\section{\bf Proof of Proposition \ref{rational}}
\label{theend}

The canonical actions of the Galois group ${\mathscr G}(\C /\Q )$ on 
$H_r\tth X\otimes_{{}_\Q}\C$ and $H^{\tth i}\tth X\otimes_{{}_\Q}\C$ 
can be transported respectively on $H_r(X,\C )$ and $H^{\tth i}(X,\C )$ 
by the canonical isomorphisms of complex vector spaces 
$$
\mu_{{}_h} : H_r\tth X\otimes_{{}_\Q}\C\to H_r(X,\C )\ \ \ \ \ \ \ \ \
\mu_{{}_c} : H^{\tth i}\tth X\otimes_{{}_\Q}\C\to H^{\tth i}(X,\C )
$$
supplied by the corresponding universal coefficients theorems.
More precisely, given $\varphi\in{\mathscr G}(\C /\Q )$
we define $\varphi^{\tth i}:H^{\tth i}(X,\C )\to H^{\tth i}(X,\C )$
by setting
$$
\varphi^{\tth i}\tth :=\tth
\mu_{{}_c}\circ (\tth id\otimes\varphi\tth )\circ\mu_{{}_c}^{-1}
$$
We will follow the usual convention to denote $\varphi^{\tth i}(\xi )$
by ${\xi}^{\tth\varphi}.$ Similar definition and notational convention
are also assumed for the homology spaces.

\smallskip

Now, the following diagram is clearly commutative
$$
\diagram
H^{\tth i}\tth X\otimes_{{}_\Q}\C\dto_{\mu_{{}_c}}
\rrto^{PD_{{}_\Q}\otimes id_{{}_\C}}&& 
H_r\tth X\otimes_{{}_\Q}\C\dto^{\mu_{{}_h}}\cr
H^{\tth i}(X,\C )\rrto_{PD_{{}_\C}}&& H_r(X,\C )
\enddiagram
$$
and it is easily seen that this implies
\begin{equation}\label{olddream}
PD_{{}_\C}\circ\varphi^{\tth i} \tth =\tth\varphi_r \circ PD_{{}_\C}
\end{equation}
i.e. the actions of ${\mathscr G}(\C /\Q )$ on $H_r(X,\C )$ and 
$H^{\tth i}(X,\C )$ are compatible with the Poincar\'e duality.

\smallskip

Take now $\xi\in S^{\tth p,i},$ and assume that $\xi$ is the 
Poincar\'e dual of $[\Gamma ]\tth ,$
where any singular simplex of $\Gamma$ satisfies (\ref{notraverso}).
Then, for every $\varphi\in{\mathscr G}(\C /\Q )$
$$
{\xi}^{\tth\varphi}\tth =\tth 
PD\bigl(\tth [\Gamma ]^{\tth\varphi}\tth\bigr)
$$
because of (\ref{olddream}). Moreover, the commutative diagram
$$
\diagram
{\mathscr Z}_r\otimes_{{}_\Q}\C\ar@{->>}[r]^-{can.}\dto_{id\otimes\varphi} &
H_r\tth X\otimes_{{}_\Q}\C\dto_{id\otimes\varphi}\rto^-{\mu_{{}_h}} &
H_r(X,\C )\dto^{\varphi_r} \cr
{\mathscr Z}_r\otimes_{{}_\Q}\C\ar@{->>}[r]_-{can.} &
H_r\tth X\otimes_{{}_\Q}\C\rto_-{\mu_{{}_h}} & H_r(X,\C ) 
\enddiagram
$$
clearly implies $[\Gamma ]^{\tth\varphi}\tth =\tth\bigl[\tth
\Gamma^{\tth\varphi}\tth\bigr]\tth .$ The cycle $\Gamma$ can be written 
as $\Gamma =\sum_h\tth w_h\tth\Gamma_h$ where $\Gamma_h\in{\mathscr Z}_r$ 
and $w_h\in\C$ for any $h\tth .$ Note that we can use the same finite set 
of singular simplexes to represent any $\Gamma_h\tth ,$ namely  
$$
\Gamma_h\tth =\tth a_{h{\sss 1}}\tth\sigma_{\sss 1}+\ldots +
a_{hs}\tth\sigma_s \ \ \ \ \ \ \ \ \ \hbox{where}\ \ \ \ a_{hj}\in\Q
\ \ \ \ \hbox{for any}\ \ \ h,j
$$
Hence
$$
\Gamma\tth =\tth\sum_{j=1}^s\tth\bigl(\tth 
a_{{\sss 1}j}\tth w_{\sss 1}+\ldots + a_{tj}\tth w_t\tth\bigr)
\tth\sigma_j
$$
Therefore our assumption on $\Gamma$ can be restated as 
follows\tth : $\sigma_j$ satisfies (\ref{notraverso}) whenever \ 
$a_{{\sss 1}j}\tth w_{\sss 1}+\ldots + a_{tj}\tth w_t\neq 0\tth .$ Finally
$$
\Gamma^{\tth\varphi}\tth =\tth\sum_h\tth\varphi (w_h)\tth\Gamma_h\tth =\tth 
\sum_{j=1}^s\tth\varphi\bigl(\tth 
a_{{\sss 1}j}\tth w_{\sss 1}+\ldots + a_{tj}\tth w_t\tth\bigr)
\tth\sigma_j
$$
which shows that also all the singular simplexes of
$\Gamma^{\tth\varphi}$ satisfy (\ref{notraverso}), 
hence ${\xi}^{\tth\varphi}\in S^{\tth p,i}.$ Then, for every
$\Q$-automorphism $\varphi$ of $\C$ we have
$(\tth S^{\tth p,i}\tth )^{\tth\varphi}\subseteq S^{\tth p,i},$
and the rationality of $S^{\tth p,i}$ follows by Galois descente.

\bigskip

\noindent{\bf Acknowledgements:} 
I wish to thank first Claire Voisin for her advice and warm 
encouragement. I owe to Elisabetta Fortuna a suggestion for
the proof of Lemma \ref{finitezza} and several explanations
concerning the semi-algebraic geometry. I thank her heartly for all this.
I would also thank Fabrizio Catanese for some useful conversations,
and Valentina Beorchia and Francesco Zucconi for their constant help
and encouragement.

\end{document}